\documentclass[11pt,table]{amsart}

\usepackage[margin=1.1in]{geometry}
\usepackage{amscd,amsmath,amsxtra,amsthm,amssymb,stmaryrd,xr,mathrsfs,mathtools,enumerate,commath,comment,orcidlink}
\usepackage{enumitem}
\usepackage{tikz}
\usetikzlibrary{arrows.meta,positioning}
\usetikzlibrary{calc}
\usepackage{tikz-3dplot}
\usepackage{stmaryrd}
\usepackage{multirow}
\usepackage{xcolor}
\usepackage{svg}
\usepackage{graphics}
\usepackage{longtable}
\usepackage{pdflscape}
\usepackage{booktabs}
\usepackage{hyperref}

\definecolor{vegasgold}{rgb}{0.77, 0.7, 0.35}
\definecolor{darkgoldenrod}{rgb}{0.72, 0.53, 0.04}
\definecolor{gold(metallic)}{rgb}{0.83, 0.69, 0.22}

\hypersetup{
 colorlinks=true,
 linkcolor=darkgoldenrod,
 filecolor=brown,
 urlcolor=gold(metallic),
 citecolor=darkgoldenrod,
}

\usepackage[all,cmtip]{xy}

\DeclareFontFamily{U}{wncy}{}
\DeclareFontShape{U}{wncy}{m}{n}{<->wncyr10}{}
\DeclareSymbolFont{mcy}{U}{wncy}{m}{n}
\DeclareMathSymbol{\Sh}{\mathord}{mcy}{"58}

\usepackage[T2A,T1]{fontenc}

\newtheorem{theorem}{Theorem}[section]
\newtheorem{lemma}[theorem]{Lemma}

\newtheorem*{theorem*}{Theorem}
\newtheorem*{ass*}{Assumption}
\newtheorem{definition}[theorem]{Definition}
\newtheorem{corollary}[theorem]{Corollary}

\newtheorem{conjecture}[theorem]{Conjecture}
\newtheorem{proposition}[theorem]{Proposition}

\newcommand{\op}[1]{\operatorname{#1}}

\newcommand{\Gal}{\op{Gal}}

\newcommand{\cE}{\mathcal{E}}

\newcommand{\ord}{\mathrm{ord}}

\newcommand{\Z}{\mathbb{Z}}
\newcommand{\Q}{\mathbb{Q}}
\newcommand{\F}{\mathbb{F}}

\newcommand{\D}{\mathfrak{D}}

\newcommand{\cS}{\mathcal{S}}

\newcommand{\Hom}{\mathrm{Hom}}

\DeclareMathSymbol{\sha}{\mathord}{mcy}{"58}

\makeatletter
\newcommand{\mylabel}[2]{#2\def\@currentlabel{#2}\label{#1}}
\makeatother

\numberwithin{equation}{section}

\newcommand{\Disc}{\operatorname{Disc}}

\newcommand{\Aut}{\operatorname{Aut}}

\newcommand{\ab}{\mathrm{ab}}

\newcommand{\cM}{\mathcal{M}}

\newcommand{\GQ}{G_{\Q}}

\newcommand{\Ind}{\operatorname{ind}}

\newcommand{\Epitop}{\operatorname{Epi}_{\mathrm{cont}}}

\title[Nilpotent Galois extensions with prescribed genus number]
{Counting nilpotent Galois extensions with prescribed genus number}

\author{Anup B. Dixit}
\address{
Department of Mathematics,
Institute of Mathematical Sciences (HBNI),
CIT Campus, IV Cross Road,
Chennai, India 600113
}
\email{anupdixit@imsc.res.in}

\author{Anwesh Ray}
\address{
Chennai Mathematical Institute,
H1, SIPCOT IT Park,
Kelambakkam, Siruseri,
Tamil Nadu 603103, India
}
\email{anwesh@cmi.ac.in}

\begin{document}
\begin{abstract}
Let $G$ be a nontrivial finite nilpotent group. In this paper, we study the distribution of genus numbers of Galois extensions of $\Q$ with Galois group isomorphic to $G$, ordered by absolute discriminant. We show that the proportion of fields with any prescribed genus number tends to zero, and obtain a quantitative bound for this proportion. Moreover, if every element of order $\ell$ in $G$ is central, where $\ell$ is the smallest prime divisor of $|G|$, we determine the orders of magnitude of the average genus number and its higher moments.
\end{abstract}



\maketitle

\section{Introduction}\label{sec:intro}

For a
number field $K/\Q$, its \emph{absolute genus field} $K^*$ is the
largest extension of $K$ obtained by adjoining an abelian extension
of $\Q$, such that the resulting extension is contained in the narrow Hilbert class field of $K$. The \emph{genus number} of $K$ is $g_K=[K^*:K]$. One can ask how the genus number varies as $K$ varies in a family of number fields of fixed degree. In this paper, we consider this question for Galois extensions of $\Q$ with a fixed non-trivial finite nilpotent Galois group $G$. More precisely, how often does the genus number attain a prescribed value in this family and how large is it on average? We also study higher moments of genus numbers in these families.\\

Before addressing these questions, let us first review counting results on the number of Galois extensions of $\Q$ with a prescribed Galois group. Shafarevich showed that every finite solvable group occurs as Galois group over $\Q$ (see
\cite[Theorem~9.6.1]{NSW}). In particular, if $G$ is nilpotent, there exists at least one Galois extension $K$ with $\op{Gal}(K/\Q)$ isomorphic to $G$. One can ask how frequently they occur when ordered by their absolute discriminant. Let $N_G(X)$ be the number of such extensions with absolute discriminant at most $X$. The strong form of Malle's conjecture \cite{Malle2002,Malle2004} predicts that
\[
 N_G(X)\sim c_G X^{a(G)}(\log X)^{b(G,\Q)-1}
\]
for a positive constant $c_G$. Here
$a(G)=\ell/((\ell-1)|G|)$, where $\ell$ is the smallest prime
divisor of $|G|$, and $b(G,\Q)$ is the number of conjugacy classes
of subgroups of order $\ell$ in $G$. The weak form predicts that for every $\varepsilon>0$,
\[
 X^{a(G)}\ll_G N_G(X)\ll_{G,\varepsilon}X^{a(G)+\varepsilon}.
\]
The weak form was proved for nilpotent groups by Kl\"uners-Malle \cite{KM}. More recently,  Koymans--Pagano constructed a parametrization of nilpotent extensions in which the coordinates record ramification. Using this parametrization, they proved the strong form whenever every  element of order $\ell$ in $G$ is central; see \cite[Theorems~1.1, 1.3 and 1.4]{KP}. In this case, $b(G,\Q)$ is the number of elements
of order $\ell$ divided by $\ell-1$. For an arbitrary finite nilpotent group, their result 
\cite[Theorem~1.3, p.~312]{KP} gives 
\[
 N_G(X)\ll_G
 X^{a(G)}(\log X)^{|I(G)|/(\ell-1)-1},
\]
where $I(G)$ is the set of elements of order $\ell$ in $G$.\\

\par We now recall some results on genus numbers in families of number fields of fixed degree. McGown-Tucker \cite{McGownTucker} studied cubic fields ordered
by absolute discriminant, treating the totally real and complex
families separately. They proved that approximately $96.23\%$
of fields in either family have genus number one. More generally, the possible genus numbers are powers of $3$ and each such power occurs for a positive proportion of cubic fields. They also showed that the average genus number tends to approximately $1.07854$, in both signatures. However, the behaviour is different if one restricts to cyclic cubic fields. In this family, the proportion with genus number one tends to zero, whereas the average genus number tends to infinity. Note that cyclic cubic fields have density zero among all cubic fields. Recently, Yamada \cite{Yamada2026} refined the cubic results of McGown-Tucker by obtaining secondary terms and sharper error estimates.\\

\par
An analogous result for quintic fields was obtained by McGown-Thorne-Tucker \cite{McGownThorneTucker}. When these fields are ordered by absolute discriminant, approximately $99.9935\%$ have genus number one, with the same proportion in each signature. Nevertheless, every fixed power of $5$ occurs as the genus number of a positive proportion of quintic fields. The limiting average is approximately $1.00026$, again independently of the signature. These results also hold on restricting to quintic fields whose Galois closures have group $S_5$.\\

\par
Kim \cite{Kim2020} studied cyclic and dihedral fields
of prime degree $\ell\geq5$. His results are unconditional for cyclic fields, while the dihedral case requires additional information about class groups. 
Bose-McGown-Panpaliya-Welling-Williams \cite{BoseEtAl}
considered pure fields of odd prime degree, namely fields
of the form $\Q(\sqrt[\ell]{a})$, ordered by absolute
discriminant. Their Theorem~1.1 shows that the proportion
with genus number one is asymptotic to a positive constant
divided by $\log X$. Their Theorem~1.2 gives an asymptotic formula for the average genus number, with order of growth $(\log X)^{\ell-1}$ and a positive leading constant.  For abelian extensions, Frei--Loughran--Newton prove a general
theorem over an arbitrary number field \cite{FLN}. They obtain an
asymptotic formula for the average genus number and show that each
fixed genus value has density zero. Their extensions on the other hand are ordered
by conductor. The first author and Pasupulati \cite{DP} studied further nonabelian
families, including degree-$3q$ fields whose Galois
closures have group $S_3\times \Z/q\Z$, where $q\neq3$ is
prime, and quartic fields whose Galois closures have
group $D_4$. Their Theorems~1.1--1.3 show that prescribed admissible genus numbers occur with positive proportions. They also obtain
asymptotic formulas for averages and higher moments in
the $S_3\times \Z/q\Z$ families. These different behaviours led them to formulate the following conjecture.

\begin{conjecture}[Dixit--Pasupulati]
\label{conj:DP}
Let $G$ be a nontrivial finite solvable group. Let $N_G(X)$
denote the number of Galois extensions $L/\Q$, counted up
to $\Q$-isomorphism, such that $\Gal(L/\Q)\simeq G$ and
$|\Disc(L)|\leq X$. For a positive integer $t$, let
$N_{G,t}(X)$ count those extensions with absolute genus
number $g_L=t$. Then, for every fixed $t\geq1$,
\[
 \lim_{X\to\infty}\frac{N_{G,t}(X)}{N_G(X)}=0.
\]
\end{conjecture}

\par
The above formulation is equivalent to
\cite[Conjecture~1.4]{DP}, which is stated in terms of towers of abelian extensions. To see this in the regular representation, note that a Galois extension $L/\Q$ admits a tower
\[
 \Q=K_0\subset K_1\subset\cdots\subset K_r=L
\]
with each $K_i/K_{i-1}$ abelian if and only if
$\Gal(L/\Q)$ is solvable. This follows from the Galois correspondence. The restriction to the regular representation is relevant here. The positive-density
results for quartic $D_4$-fields concern non-Galois fields of degree four, whereas a Galois extension with group $D_4$ has degree eight.

\subsection{Main results}

We now state the main results. Let $G$ be a  nontrivial finite nilpotent group, and let $\ell$ be the smallest prime divisor of its order. Recall that $N_G(X)$ counts the Galois extensions $L/\Q$ with $\Gal(L/\Q)\simeq G$ and absolute discriminant at most $X$. Let $N_{G,t}(X)$ count the number of these extensions with genus number $t$. Let $d_\ell(G)$ be the largest integer $d$ for which
$G$ has a quotient isomorphic to $(\Z/\ell\Z)^d$, and let
$\nu_\ell(t)$ denote the exponent of $\ell$ in the prime factorization
of $t$. Set $R_G(t)=d_\ell(G)+\nu_\ell(t)$. We prove the following.

\begin{theorem}
\label{thm:main}
For every fixed positive integer $t$ and every $X\geq3$,
\begin{equation}\label{eq:main-intro}
 N_{G,t}(X)
 \ll_{G,t}
 \frac{X^{a(G)}}{\log X}
 (1+\log\log X)^{R_G(t)-1}.
\end{equation}
\end{theorem}
By the lower bound of Kl\"uners-Malle, we have $N_G(X) \gg_G X^{a(G)}$. Thus, we immediately deduce the following consequence.

\begin{corollary}\label{cor:density}
For every nontrivial finite nilpotent group $G$ and every fixed
positive integer $t$, the proportion of Galois $G$-extensions of
$\Q$ with genus number $t$, ordered by absolute discriminant,
tends to zero.
\end{corollary}

This proves Conjecture \ref{conj:DP} specialized to finite nilpotent groups. The proofs of these results are given in Section~\ref{sec:counting}.\\

\par Our next result concerns moments of genus numbers. For a real number ${r}\geq0$, define
\[
 M_{G,{r}}(X)
 :=
 \sum_{\substack{
 L/\Q\ \mathrm{Galois}\\
 \Gal(L/\Q)\simeq G\\
 |\Disc(L)|\leq X
 }}
 g_L^{r},
\]
where the sum runs over $\Q$-isomorphism classes of extensions. Let $I(G)$ denote the set of elements of order $\ell$ in $G$. If every element in $I(G)$ is central, then $b(G,\Q)=|I(G)|/(\ell-1)$. Under this assumption, we determine the order of growth of $M_{G,{r}}(X)$.

\begin{theorem}\label{thm:moments-intro}
Suppose that every element of order $\ell$ in $G$ is central.
For every fixed real ${r}\geq0$,
\begin{equation}\label{eq:moment-central-intro}
 M_{G,{r}}(X)
 \asymp_{G,{r}}
 X^{a(G)}(\log X)^{\ell^{r} b(G,\Q)-1}.
\end{equation}
\end{theorem}

Note that $M_{G,0}(X) = N_G(X)$. Applying the theorem with ${r}=0$ and dividing the estimates, we obtain for any ${r}\geq 0$,
\[
 \frac{M_{G,{r}}(X)}{N_G(X)}
 \asymp_{G,{r}}(\log X)^{(\ell^{r}-1)b(G,\Q)}.
\]
In particular, the average genus number has order
$(\log X)^{(\ell-1)b(G,\Q)}$. For illustration, consider $G=\Z/\ell\Z$. In this case, $d_\ell(G)=b(G,\Q)=1$, and hence for every fixed integer $k\geq 0$,
\[
 N_{\Z/\ell\Z,\ell^k}(X)
 \ll_{\ell,k}\frac{X^{1/(\ell-1)}}{\log X}
 (1+\log\log X)^k,\]
 and
\[ \frac{M_{\Z/\ell\Z,{r}}(X)}{N_{\Z/\ell\Z}(X)}
 \asymp_{\ell,{r}}(\log X)^{\ell^{r}-1}.
\]
The quaternion group $Q_8$ provides a nonabelian example to which the moment theorem applies. On the other hand, the dihedral group $D_4$ of order eight
satisfies the hypothesis of Theorem~\ref{thm:main}, but not the
centrality assumption in Theorem~\ref{thm:moments-intro}.

\subsection{Methodology}

The key ingredient in our approach is the relation between the genus number and the tamely ramified primes whose inertia groups have order $\ell$. For each such prime, there is an associated cyclic
extension of degree $\ell$ inside a cyclotomic field. Adjoining
this extension to $L$ introduces no new ramification at finite places, so it gives a subextension of the genus field. The
cyclotomic subfields coming from distinct primes are linearly
disjoint. Consequently, these primes force a growing power of $\ell$ to divide
the genus number. In particular, fixing the genus number bounds the number of such primes. They are also the ramified primes with the
smallest possible contribution to the discriminant. We use Koymans-Pagano parameterization to count extensions subject to this restriction. The solvable tuples parametrize epimorphisms bijectively and each field corresponds to exactly $|\Aut(G)|$ epimorphisms. 
\subsection{Organization}

The paper is organized as follows. In Section~2, we collect the preliminary results on discriminants
and genus fields. We construct cyclotomic extensions inside the genus field
and use their linear disjointness to obtain lower bounds for
the genus number. Section~\ref{sec:KP-local} recalls the Koymans-Pagano coordinates,
first for groups of prime-power order and then for a general
nilpotent group by combining its Sylow factors. We explain how
these coordinates capture tame inertia and the discriminant.
We then prove the divisibility statement which bounds the number
of primes with inertia of order $\ell$ when the genus number is
fixed. The section concludes with the estimates for integers with a
bounded number of prime divisors that are needed in the counting
argument. In Section~\ref{sec:counting}, we prove Theorem~\ref{thm:main} and Corollary \ref{cor:density}. Section~\ref{sec:weighted} is devoted to the weighted counting results.
We first bound the genus number in terms of ramification indices and obtain a weighted upper bound for every nilpotent group. The argument treats the finitely many exceptional primes separately and applies Selberg-Delange estimates to the remaining primes. A central-twisting construction then gives a matching lower bound under the centrality assumption. Finally, we deduce the orders of magnitude of the averages and higher moments of genus numbers.

\subsection*{Declarations}
\begin{description}
     \item[Conflict of interest] There is no conflict of interest to report.
    \item[Availability of data and materials] No data was generated or analyzed in obtaining the results in this article.
    \item[Funding] No funding was received during the preparation of this article.
    \item[Use of artificial intelligence]
The authors used ChatGPT 6.1 Sol as an exploratory tool to study related papers and to assist in correcting mistakes in the article. 
\end{description}

\section{Preliminaries}\label{sec:preliminaries}

Let a nontrivial finite group $G$ act on itself by left multiplication.
For $g\in G$, let $\Ind(g)$ denote the permutation index of $g$, namely
\begin{equation}\label{ing g defn}
 \Ind(g)=|G|-\#\{\text{orbits of }\langle g\rangle\text{ on }G\}.
\end{equation}
\noindent Suppose that \(g\in G\) has order \(m\). Then every orbit of the
cyclic subgroup \(\langle g\rangle\) has cardinality exactly \(m\).
Consequently, the regular action of \(g\) decomposes \(G\) into exactly
\(|G|/m\) cycles and it follows that
\begin{equation}\label{eq:regular-index}
        \Ind(g)
        =
        |G|-\frac{|G|}{m}
        =
        |G|\left(1-\frac1m\right).
\end{equation}

Let \(\ell\) denote the smallest prime divisor of \(|G|\).
Every nonidentity element of $G$ has order at least $\ell$, and
Cauchy's theorem guarantees the existence of an element of order $\ell$. Consequently, we find that
\begin{equation}\label{eq:minimal-regular-index}
 \min_{1\neq g\in G}\Ind(g)
 =
 |G|\left(1-\frac1\ell\right).
\end{equation}
\begin{definition}\label{def:eg}
For $g\neq1$ set
\[
 e_g:=\Ind(g)=|G|\left(1-\frac1{\ord(g)}\right).
\]
Also set
\[
 e_0:=|G|\left(1-\frac1\ell\right)
 \quad\text{and}\quad
 a(G):=e_0^{-1}=\frac{\ell}{(\ell-1)|G|}.
\]
\end{definition}
\noindent We note that if $\ord(g)>\ell$ then it follows that
$e_g>e_0$, and consequently we have that $a(G)e_g>1$.

\begin{lemma}\label{lem:tame-disc}
Let $L/K$ be a finite Galois extension of number fields of degree $n$.
Let $\mathfrak q$ be a finite prime of $K$ which is tamely ramified in
$L$, and let $e$ be the common ramification index of the primes of $L$
above $\mathfrak q$. Then
\begin{equation}\label{eq:tame-disc-general-base}
 v_{\mathfrak q}(\mathfrak D_{L/K})
 =n\left(1-\frac1e\right),
\end{equation}
where $\mathfrak D_{L/K}$ is the relative discriminant ideal.
\end{lemma}

\begin{proof}
By \cite[Chapter~IV, \S1, Proposition~4]{SerreLocalFields}, the
different exponent at every prime $w\mid\mathfrak q$ is $e-1$.
Writing $f_w$ for the residue degree at $w$ and taking the norm of
the different, we obtain
\[
 v_{\mathfrak q}(\mathfrak D_{L/K})
 =
 \sum_{w\mid\mathfrak q}f_w(e-1)
 =
 n\left(1-\frac1e\right),
\]
since $\sum_{w\mid\mathfrak q}ef_w=n$.
\end{proof}

In particular, if $L/\Q$ has Galois group $G$ and $g\in G$
generates the tame inertia group at a rational prime $q$, then
$e_q(L/\Q)=\ord(g)$. The exponent of $q$ in the discriminant is
therefore
\[
 |G|\left(1-\frac1{\ord(g)}\right)=e_g.
\]

Fix an algebraic closure $\overline\Q$ and write
$\GQ=\Gal(\overline\Q/\Q)$. For a finite group $G$, let
$\Epitop(\GQ,G)$ denote the continuous group epimorphisms
$\GQ\twoheadrightarrow G$.

\begin{lemma}\label{lem:aut-multiplicity}
Let $\mathscr F$ be any finite set of Galois extensions $L/\Q$
with $\Gal(L/\Q)\simeq G$, where membership in $\mathscr F$
depends only on the field $L$ and not on a chosen isomorphism
$\Gal(L/\Q)\simeq G$. Then the number of epimorphisms
$\psi:\GQ\twoheadrightarrow G$ whose kernel has fixed field
belonging to $\mathscr F$ is $|\Aut(G)|\,|\mathscr F|$.
\end{lemma}

\begin{proof}
For a fixed field $L$, the epimorphisms in question correspond
exactly to isomorphisms $\Gal(L/\Q)\simeq G$. There are
$|\Aut(G)|$ such isomorphisms. Summing over
$L\in\mathscr F$ proves the assertion.
\end{proof}

\begin{definition}\label{def:genus}
Let $K/\Q$ be a number field. The \emph{absolute genus field}
$K^*$ is the maximal extension of $K$ satisfying:
\begin{enumerate}[label=\rm(\roman*)]
\item $K^*/K$ is unramified at every finite place of $K$;
\item there exists a finite abelian extension $F/\Q$ such that
$K^*=KF$.
\end{enumerate}
The \emph{genus number} is $g_K=[K^*:K]$.
\end{definition}

We allow ramification at the infinite places and this convention for the
absolute genus field is used throughout. Every extension $KF/K$ satisfying
\textup{(i)} and \textup{(ii)} is abelian over $K$: restriction
identifies its Galois group with a subgroup of $\Gal(F/\Q)$.
It therefore lies in the narrow Hilbert class field of $K$.
Such extensions are closed under compositum, and their
compositum is finite because the narrow Hilbert class field
is finite. This gives a unique largest extension of the
required form. In particular, $K^*/K$ is finite.\\

We also recall Abhyankar's lemma (see \cite[Chapter 5, Lemma 3]{Narkiewicz}).

\begin{lemma}\label{lem:abhyankar}
Let $F$ be a nonarchimedean local field with residue
characteristic $p$. Let $M/F$ and $N/F$ be finite extensions
inside a common algebraic closure of $F$. Assume that $M/F$
is tamely ramified and that $e(M/F)\mid e(N/F)$.
Then $MN/N$ is unramified.
\end{lemma}

Let $q$ be a rational prime and let $d\mid q-1$. Since
\[
 \Gal(\Q(\zeta_q)/\Q)\simeq(\Z/q\Z)^\times
\]
is cyclic of order $q-1$, there is a unique subfield
$E_{q,d}\subseteq\Q(\zeta_q)$ with $[E_{q,d}:\Q]=d$.
The extension $E_{q,d}/\Q$ is cyclic, ramified only at $q$
among finite primes, and totally tamely ramified at $q$
with ramification index $d$.\\

\par Now let $K/\Q$ be Galois and suppose that $q$ is tamely
ramified in $K$. Let $e_q=e_q(K/\Q)$ be its ramification index
and set $d_q:=\gcd(e_q,q-1)$ and $E_q:=E_{q,d_q}$.

\begin{proposition}\label{prop:oneprime}
With the notation above, $KE_q/K$ is unramified at every
finite place. Hence $KE_q\subseteq K^*$ and
\[
 [KE_q:K]=[E_q:E_q\cap K]\mid g_K.
\]
In particular,
\[
 g_K\geq[E_q:E_q\cap K].
\]
\end{proposition}

\begin{proof}
It is clear that $KE_q/K$ is unramified at every prime of $K$
that does not lie above $q$. Let $v$ be the unique prime of
$E_q$ above $q$. For each prime $w\mid q$ of $K$, place the
completions $(E_q)_v$ and $K_w$ inside a common algebraic
closure of $\Q_q$. The extension $(E_q)_v/\Q_q$ is tamely
ramified of ramification index $d_q$, and by definition
\[
 d_q\mid e_q=e(K_w/\Q_q).
\]
Lemma~\ref{lem:abhyankar} therefore shows that
$K_w(E_q)_v/K_w$ is unramified. Hence $KE_q/K$ is unramified
at every prime above $q$ as well. Thus $KE_q\subseteq K^*$,
and the asserted divisibility follows from the tower law.
\end{proof}

\begin{lemma}\label{lem:cyclo-disjoint}
Let $q_1,\dots,q_r$ be distinct rational primes, and for each
$i$ let $E_i/\Q$ be a nontrivial subextension of
$\Q(\zeta_{q_i})/\Q$. Then the fields $E_1,\dots,E_r$ are
linearly disjoint over $\Q$. In particular,
\[
 [E_1\cdots E_r:\Q]=\prod_{i=1}^r[E_i:\Q].
\]
\end{lemma}

\begin{proof}
We argue inductively. Set $F=E_1\cdots E_{r-1}$.
The extension $F/\Q$ is unramified at $q_r$, since each
$E_i$ with $i<r$ is ramified only at $q_i$. On the other hand,
$E_r/\Q$ is totally ramified at $q_r$. Let
\[
 D=E_r\cap F.
\]
As a subfield of $F$, the extension $D/\Q$ is unramified at
$q_r$. As a subfield of the cyclic totally ramified extension
$E_r/\Q$, it is totally ramified at $q_r$: the inertia group
of $E_r/\Q$ is the whole Galois group, and its image in every
quotient is again the whole quotient. The only extension
which is simultaneously unramified and totally ramified at
$q_r$ is the trivial extension. Thus $D=\Q$. Hence $E_r$ and
$F$ are linearly disjoint, and induction completes the proof.
\end{proof}

We write $G^{\ab}:=G/[G,G]$ for the abelianization of $G$.

\begin{proposition}\label{prop:tame-lower}
Let $K/\Q$ be a finite Galois extension with group $G$, and
let $T$ be any finite set of rational primes tamely ramified
in $K$. Setting $E_T:=\prod_{q\in T}E_q$, where the empty
compositum is understood to be $\Q$, we have that $KE_T/K$
is unramified at every finite place and
\[
 [KE_T:K]\mid g_K.
\]
Moreover,
\begin{equation}\label{eq:tame-lower}
 g_K\geq[KE_T:K]
 =\frac{\prod_{q\in T}d_q}{[E_T\cap K:\Q]}
 \geq\frac1{|G^{\ab}|}\prod_{q\in T}d_q.
\end{equation}
\end{proposition}

\begin{proof}
By Proposition~\ref{prop:oneprime}, each $KE_q/K$ is
unramified and therefore $KE_T\subseteq K^*$. Consequently,
$[KE_T:K]\mid g_K$. After omitting the trivial factors
$E_q=\Q$, Lemma~\ref{lem:cyclo-disjoint} implies that
\[
 [E_T:\Q]=\prod_{q\in T}d_q.
\]
Therefore
\[
 [KE_T:K]=[E_T:E_T\cap K]
 =\frac{\prod_{q\in T}d_q}{[E_T\cap K:\Q]}.
\]
Finally $E_T\cap K$ is Galois and abelian over $\Q$, because
it is the intersection of two Galois extensions and is
contained in the abelian extension $E_T/\Q$. Hence
$\Gal(E_T\cap K/\Q)$ is an abelian quotient of $G$, so
\[
 [E_T\cap K:\Q]\leq |G^{\ab}|.
\]
This proves \eqref{eq:tame-lower}.
\end{proof}

\begin{definition}\label{def:dell}
For a prime $p$ and a finite group $G$, define
\[
 d_p(G):=\dim_{\F_p}\Hom(G,\F_p),
\]
where $\F_p$ is viewed as an additive group.
\end{definition}
\noindent Note that $p^{d_p(G)}$ is the order of the maximal
elementary abelian $p$-quotient of $G$. This quotient is $G/\bigl([G,G]G^p\bigr)$, where $G^p:=\langle g^p:g\in G\rangle$. An elementary abelian $p$-extension is a finite Galois
extension with Galois group isomorphic to $(\Z/p\Z)^r$ for
some $r\geq0$.

\begin{lemma}\label{lem:intersection-rank}
Let $K/\Q$ be Galois with group $G$, and let $E/\Q$ be an
elementary abelian $p$-extension. Then
\[
 [E\cap K:\Q]\leq p^{d_p(G)}.
\]
\end{lemma}

\begin{proof}
The field $D=E\cap K$ is Galois over $\Q$. Since $D\subseteq E$,
its Galois group is an elementary abelian $p$-group. Since
$D\subseteq K$, restriction identifies $\Gal(D/\Q)$ with a
quotient of $G$. If $[D:\Q]=p^r$, then $G$ therefore has a
quotient $(\Z/p\Z)^r$, so $r\leq d_p(G)$ by
Definition~\ref{def:dell}. This proves the claim.
\end{proof}

\begin{proposition}\label{prop:prime-order}
Let $K/\Q$ be Galois with group $G$. Let $T$ be a finite set
of tamely ramified primes such that $e_q(K/\Q)=\ell$ and
$q\equiv1\pmod\ell$ for every $q\in T$. Then
\begin{equation}\label{eq:prime-order-lower}
 \nu_\ell(g_K)\geq\max\{0,|T|-d_\ell(G)\}.
\end{equation}
Equivalently,
\[
 \ell^{\max\{0,|T|-d_\ell(G)\}}\mid g_K.
\]
In particular,
\[
 g_K\geq\ell^{|T|-d_\ell(G)}.
\]
\end{proposition}

\begin{proof}
For $q\in T$, one has
\[
 d_q=\gcd(\ell,q-1)=\ell.
\]
Thus $E_q/\Q$ is cyclic of degree $\ell$.
By Lemma~\ref{lem:cyclo-disjoint},
\[
 \Gal(E_T/\Q)\simeq(\Z/\ell\Z)^{|T|},
\]
and consequently $
 [E_T:\Q]=\ell^{|T|}$.
Write $[E_T\cap K:\Q]=\ell^h$. By
Lemma~\ref{lem:intersection-rank}, we have $h\leq d_\ell(G)$.
Proposition~\ref{prop:tame-lower} gives
\[
 [KE_T:K]
 =
 [E_T:E_T\cap K]
 =
 \ell^{|T|-h}
 \mid g_K.
\]
It follows that
\[
 \nu_\ell(g_K)\geq |T|-h\geq |T|-d_\ell(G).
\]
Combining this with $\nu_\ell(g_K)\geq0$ proves
\eqref{eq:prime-order-lower}.
\end{proof}

The congruence condition in Proposition~\ref{prop:prime-order}
follows automatically from the other hypotheses, since
$\ell$ is the smallest prime divisor of $|G|$. Indeed, choose
a prime of $K$ above $q$, and denote its decomposition and
inertia groups by $D_q$ and $I_q$. Since $I_q\simeq \Z/\ell\Z$,
the image of the conjugation action
\[
 D_q\longrightarrow\Aut(I_q)
\]
has order dividing both $|G|$ and $\ell-1$. These integers
are coprime, so this action is trivial. On the other hand,
the tame Frobenius relation gives
$\sigma\tau\sigma^{-1}=\tau^q$, where $\tau$ generates $I_q$
and $\sigma$ lifts arithmetic Frobenius. Thus
$\tau^{q-1}=1$, which gives $q\equiv1\pmod\ell$.
\section{Koymans--Pagano coordinates and the genus--inertia dictionary}
\label{sec:KP-local}
\label{sec:KP-power}
\label{sec:KP-nilpotent}
\label{sec:dictionary}
\label{sec:analytic}

\par We recall the Koymans--Pagano parametrization
\cite[Section~2]{KP}. This will be needed in the proof of our fixed-genus
counting results as well as our results on moments. The parameterization describes nilpotent extensions in terms of squarefree ideal coordinates which record the image of tame
inertia. Throughout, we work over $\Q$. Let $G$ be a nontrivial finite nilpotent group.
Write $G=\prod_{j=1}^cP_j$, where $P_j=G(\ell_j)$ is the Sylow
$\ell_j$-subgroup and the $\ell_j$ are the distinct prime divisors of
$|G|$. We order these primes so that
$\ell=\ell_1=\min\{p:p\mid |G|\}$.
For each $j$, let $G_{\Q}^{\mathrm{pro}-\ell_j}$ denote the maximal
pro-$\ell_j$ quotient of $G_{\Q}$.\\

By \cite[Section~2]{KP}, there is a finite set of rational
primes $S_G$, containing every prime divisor of $|G|$, for which the
following choices can be made. Set
\[
 \widetilde\Omega_{\Q}(\ell_j)
 :=\{q\notin S_G:q\equiv1\pmod{\ell_j}\}.
\]
For every $q\in\widetilde\Omega_{\Q}(\ell_j)$, fix a character
\[
 \chi_{j,q}:G_{\Q}\longrightarrow\F_{\ell_j}
\]
which is ramified at $q$ and unramified at every finite prime outside
$S_G\cup\{q\}$. Also fix an $\F_{\ell_j}$-basis
$\chi_{j,1},\dots,\chi_{j,t_j}$ of the finite-dimensional space of
characters unramified at every finite prime outside $S_G$, so that
\[
 \{\chi_{j,1},\dots,\chi_{j,t_j}\}
 \cup
 \{\chi_{j,q}:q\in\widetilde\Omega_{\Q}(\ell_j)\}
\]
is a basis of $H^1(G_{\Q},\F_{\ell_j})$. Here $\F_{\ell_j}$ is
regarded as a trivial Galois module, so these cohomology classes are
continuous homomorphisms to its additive group. When $j=1$, we write
$\chi_q:=\chi_{1,q}$.\\

For $q\notin S_G$, let $I_q(\ell_j)$ denote the image of inertia in
$G_{\Q}^{\mathrm{pro}-\ell_j}$, after fixing compatible local
embeddings. If $q\in\widetilde\Omega_{\Q}(\ell_j)$, choose a
topological generator $\sigma_{j,q}$ of $I_q(\ell_j)$ with
$\chi_{j,q}(\sigma_{j,q})=1$. Otherwise $I_q(\ell_j)$ is trivial,
and we set $\sigma_{j,q}=1$. The latter assertion concerns the image
of inertia in the pro-$\ell_j$ Galois quotient. The Frobenius
relation with a tame inertia generator forces this image to be trivial when
$q\not\equiv1\pmod{\ell_j}$.\\

\par Let $\cS'_j$ denote the set of squarefree positive integers
supported entirely on $\widetilde\Omega_{\Q}(\ell_j)$, and set
\[
 \cS_j:=\{0,1\}^{[t_j]}\times\cS'_j\quad \text{where}\quad [t_j]:=\{1,\dots,t_j\}.
\]
Two
elements $w=(T,b)$ and $w'=(T',b')$ are called \emph{coprime} if
$b$ and $b'$ are relatively prime and there is no $i\in[t_j]$ for
which $T_i=T'_i=1$. We write
$\operatorname{Prim}(\cS_j^{P_j\setminus\{1\}})$ for the set of
tuples $(w_{h,j})_{h\in P_j\setminus\{1\}}$ whose coordinates are
pairwise coprime. 

\par Projection onto the Sylow factors gives rise to the following identification:
\[
 \Epitop(G_{\Q},G)
 \simeq\prod_{j=1}^c\Epitop(G_{\Q}^{\mathrm{pro}-\ell_j},P_j);
\]
see \cite[Proposition~2.11]{KP}. We use the common exceptional set
$S_G\cup\{\infty\}$ in every Sylow construction. Enlarging the
original exceptional sets only enlarges the finite exceptional
character spaces.\\

Koymans--Pagano bundle the Sylow coordinates into a tuple
$(v_g(1),v_g(2))_{g\in G\setminus\{1\}}$ as follows. Set
$w_{1,j}:=(0,1)$. For $q\notin S_G$, let $g_j(q)=h$ if
$q\mid w_{h,j}(2)$, and let $g_j(q)=1$ if $q$ divides none of the
$w_{h,j}(2)$. For $g=(g_1,\dots,g_c)\ne1$, set
\[
 v_g(1):=(w_{g_1,1}(1),\dots,w_{g_c,c}(1))\quad \text{and}\quad
 v_g(2):=
 \prod_{\substack{q\notin S_G\\(g_1(q),\dots,g_c(q))=g}}q.
\]
The $v_g(2)$ are squarefree, pairwise coprime, and supported outside
$S_G$. The first components take values in fixed finite sets. Let $\cS'$ be the set of squarefree positive integers supported outside
$S_G$, and set \[\cS:=\{0,1\}^{[t_1+\cdots+t_c]}\times\cS'.\]
We denote the image of the product of the Sylow primitive spaces
under this bundling by
$\operatorname{Prim}(\cS^{G\setminus\{1\}})$.
Fix the construction choices in \cite[Section~2]{KP}. The
\emph{solvable primitive tuples}, denoted by
$\operatorname{Prim}(\cS^{G\setminus\{1\}})^{\mathrm{solv}}$ are those that give rise to epimorphisms onto $G$.

\begin{theorem}[Koymans--Pagano parametrization]\label{thm:KP-l-param}
There is a bijection
\begin{equation}\label{eq:KP-l-bijection}
 P_G:
 \operatorname{Prim}(\cS^{G\setminus\{1\}})^{\mathrm{solv}}
 \stackrel{\sim}{\longrightarrow}
 \Epitop(G_{\Q},G).
\end{equation}
\end{theorem}

\begin{proof}
This follows from \cite[Propositions~2.8 and~2.11]{KP} and the
bundling described above.
\end{proof}

There is also a restriction on which primes may occur in a given
coordinate. If $q\mid v_g(2)$ and $\ell_j\mid\ord(g)$, then the
$\ell_j$-part of inertia is nontrivial. Hence
$q\in\widetilde\Omega_{\Q}(\ell_j)$ and therefore
$q\equiv1\pmod{\ell_j}$. In particular, whenever $\ord(g)=\ell$ and
$q\mid v_g(2)$, one has $q\equiv1\pmod\ell$.

\begin{proposition}\label{prop:bundled-inertia}
Let $v$ be a solvable bundled tuple corresponding to a $G$-extension
$L/\Q$, and let $q\notin S_G$. If $q\mid v_g(2)$, then a tame inertia
generator at $q$ maps to $g$. Consequently, the ramification index
satisfies $e_q(L/\Q)=\ord(g)$. If $q$ divides none of the $v_g(2)$,
then $q$ is unramified in $L$.
\end{proposition}

\begin{proof}
This is \cite[Proposition~2.12]{KP}.
\end{proof}

For a positive integer $d$ and a finite set $S$ of rational primes,
write
\[
 \operatorname{free}_S(d):=\prod_{q\notin S}q^{v_q(d)}
\]
for the largest divisor of $d$ supported outside $S$.

\begin{proposition}
\label{prop:KP-disc-general}
For a solvable bundled tuple corresponding to $L/\Q$,
\begin{equation}\label{eq:KP-disc-general}
 \operatorname{free}_{S_G}|\Disc(L)|
 =
 \prod_{g\in G\setminus\{1\}}v_g(2)^{e_g},
\end{equation}
where we recall that $e_g=|G|\left(1-\frac1{\ord(g)}\right)$. Consequently, if $|\Disc(L)|\leq X$, then
\begin{equation}\label{eq:disc-upper-monomial}
 \prod_{g\neq1}v_g(2)^{e_g}\leq X.
\end{equation}
\end{proposition}

\begin{proof}
The result follows from \cite[Proposition~2.13]{KP}.
\end{proof}

Now let
\[I(G):=\{g\in G\setminus\{1\}:\ord(g)=\ell\}.\]
By Section~\ref{sec:preliminaries}, the elements of $I(G)$ are exactly
the nonidentity elements of minimal index, and
\[
 e_h=e_0=a(G)^{-1}\quad(h\in I(G)),
 \qquad a(G)e_g>1\quad(\ord(g)>\ell).
\]

\begin{proposition}
\label{prop:dictionary}
Let $L/\Q$ be a nilpotent Galois extension with group $G$, represented
by a solvable bundled Koymans--Pagano tuple.  Let $q\notin S_G$ and
suppose that $q\mid v_g(2)$.  Then
\begin{enumerate}[label=\rm(\roman*)]
\item the tame inertia group at $q$ has order $\ord(g)$;
\item the exponent of $q$ in $\Disc(L)$ is
\[
        |G|\left(1-\ord(g)^{-1}\right);
\]
\item the degree of the tame cyclotomic extension associated to $q$
in Ishida's construction \cite{Ishida} is
\[
        \gcd(\ord(g),q-1).
\]
\end{enumerate}
In particular, if $g\in I(G)$, then this degree is $\ell$.
\end{proposition}

\begin{proof}
Parts~(i) and~(ii) follow from
Propositions~\ref{prop:bundled-inertia} and
\ref{prop:KP-disc-general}.  For a tamely ramified prime $q$,
the cyclotomic construction preceding
Proposition~\ref{prop:oneprime} associates the unique subfield of
$\Q(\zeta_q)$ of degree $\gcd(e_q(L/\Q),q-1)$.
By~(i), one has $e_q(L/\Q)=\ord(g)$,
which proves~(iii). Finally, if $g\in I(G)$, then $\ord(g)=\ell$.  Moreover,
$q\mid v_g(2)$ implies that $q$ belongs to the clean set supporting
the $\ell$-coordinates, and hence $q\equiv1\pmod\ell$.  Therefore
\[
        \gcd(\ord(g),q-1)
        =\gcd(\ell,q-1)
        =\ell.
\]
\end{proof}

The number of primes appearing in the minimal coordinates admits an
intrinsic description.  For a solvable tuple $v$, set
$M(v):=\prod_{h\in I(G)}v_h(2)$.  Since the integers $v_h(2)$ are
squarefree and pairwise coprime, $M(v)$ is squarefree and
$\omega(M(v))=\sum_{h\in I(G)}\omega(v_h(2))$, where $\omega(n)$ denotes the number of distinct prime divisors of a squarefree integer $n$. By Proposition~\ref{prop:bundled-inertia}, this quantity depends only
on the extension and not on the chosen coordinates.  We therefore
define
\[
 \omega_{\min,S_G}(L)
 :=
 \#\{q\notin S_G:
 q\text{ ramifies in }L,\ e_q(L/\Q)=\ell\}.
\]
Thus $\omega_{\min,S_G}(L)=\omega(M(v))$.  In words,
$\omega_{\min,S_G}(L)$ counts the clean ramified primes whose inertia
group has the smallest possible nontrivial order.  These are exactly
the primes which occur in the minimal-index coordinates.

\begin{theorem}
\label{thm:genus-minimal}
Let $L/\Q$ be a nilpotent Galois extension with group $G$, and set
$r:=\omega_{\min,S_G}(L)$.  Then
\begin{equation}\label{eq:genus-divisibility}
 \nu_\ell(g_L)\geq
 \max\{0,r-d_\ell(G)\}.
\end{equation}
In particular, if $r>d_\ell(G)$, then
$\ell^{\,r-d_\ell(G)}\mid g_L$, and hence
$g_L\geq\ell^{\,r-d_\ell(G)}$.
\end{theorem}

\begin{proof}
Let $T$ denote the set of primes counted by
$\omega_{\min,S_G}(L)$, so that $|T|=r$. These primes are
tamely ramified in $L$, and Proposition~\ref{prop:dictionary}
gives
\[
 e_q(L/\Q)=\ell
 \quad\text{and}\quad
 q\equiv1\pmod\ell
 \qquad(q\in T).
\]
Proposition~\ref{prop:prime-order} therefore yields
\[
 \nu_\ell(g_L)\geq\max\{0,r-d_\ell(G)\}.
\]
The remaining assertions follow immediately.
\end{proof}

Thus fixing the genus number imposes a uniform bound on the number of
prime divisors which may occur among all the minimal coordinates.

\begin{corollary}
\label{cor:omega-fixed}
For $t\geq1$, set $R_G(t):=d_\ell(G)+\nu_\ell(t)$.  If $g_L=t$, then
$\omega_{\min,S_G}(L)\leq R_G(t)$.  Equivalently, every
Koymans--Pagano tuple representing $L$ satisfies
\[
 \omega\!\left(\prod_{h\in I(G)}v_h(2)\right)
 \leq R_G(t).
\]
\end{corollary}

\begin{proof}
Set $r=\omega_{\min,S_G}(L)$.  By
Theorem~\ref{thm:genus-minimal},
$\max\{0,r-d_\ell(G)\}\leq\nu_\ell(g_L)=\nu_\ell(t)$.
Thus $r\leq d_\ell(G)+\nu_\ell(t)=R_G(t)$.
\end{proof}

Notice that $R_G(t)$ is always positive.  Indeed, the Sylow
$\ell$-subgroup of $G$ is nontrivial and, as every nontrivial finite
$\ell$-group has a quotient of order $\ell$, projection to this Sylow
factor followed by such a quotient gives a nonzero homomorphism
$G\to\F_\ell$.  Thus $d_\ell(G)\geq1$, and consequently
$R_G(t)\geq1$ for every $t\geq1$.\\

For an integer $r\geq0$ and $Y\geq1$, let
$A_r(Y):=\#\{n\leq Y:\mu^2(n)=1,\ \omega(n)=r\}$.  Here $\mu$ is the
Möbius function, so the condition $\mu^2(n)=1$ is equivalent to $n$
being squarefree, while $\omega(n)$ denotes the number of distinct
prime divisors of $n$. For our application, we need the following upper bound.

\begin{lemma}
\label{lem:colored-landau}
Fix integers $M\geq1$ and $R\geq1$, and define
\[
 B_{M,R}(Y)
 :=
 \sum_{\substack{n\leq Y\\
                  \mu^2(n)=1\\
                  \omega(n)\leq R}}
 M^{\omega(n)}.
\]
Then
\begin{equation}\label{eq:colored-landau}
 B_{M,R}(Y)
 \ll_{M,R}
 \frac{Y}{\log Y}
 (1+\log\log Y)^{R-1}
\end{equation}
for $Y\geq 3$. For all $Y\geq1$ one also has the bound
$B_{M,R}(Y)\ll_{M,R}Y$.
\end{lemma}

\begin{proof}
We first prove that, for every fixed integer $r\geq1$,
\begin{equation}\label{bound}
 A_r(Y)
 \ll_r
 \frac{Y}{\log Y}(1+\log\log Y)^{r-1},
\end{equation}
for $Y\geq 3$. For $r=1$, the integers counted by $A_1(Y)$ are precisely
the primes, so this follows from Chebyshev's estimate
$\pi(Y)\ll Y/\log Y$. Suppose that $r\geq2$ and that the assertion holds for $r-1$.
Every squarefree integer $n\leq Y$ with exactly $r$ prime
divisors has a smallest prime divisor $p\leq Y^{1/r}$.
Writing $n=pm$, the integer $m$ is squarefree and has
exactly $r-1$ prime divisors. Consequently,
\[
 A_r(Y)
 \leq
 \sum_{p\leq Y^{1/r}}A_{r-1}(Y/p),
\]
where the sum is over primes. This is an upper bound
because the sum also allows integers $m$ divisible by $p$
or having a prime divisor smaller than $p$. Assume first that $Y\geq9$. For every prime in this sum,
\[
 Y/p\geq Y^{1-1/r}\geq\sqrt Y\geq3\quad \text{and thus}\quad 
 \log(Y/p)\geq\frac12\log Y.
\]
The induction hypothesis therefore gives
\[
 A_r(Y)
 \ll_r
 \frac{Y}{\log Y}
 (1+\log\log Y)^{r-2}
 \sum_{p\leq Y^{1/r}}\frac1p.
\]
Chebyshev's estimate and partial summation imply
\[
 \sum_{p\leq Y}\frac1p
 =
 \frac{\pi(Y)}{Y}
 +\int_2^Y\frac{\pi(t)}{t^2}\,dt
 \ll 1+\log\log Y.
\]
Since the prime sum over $p\leq Y^{1/r}$ is no larger,
we obtain the required bound for $A_r(Y)$. Enlarging
the implied constant covers $3\leq Y<9$ and completes
the induction.\\

We now group the terms defining $B_{M,R}(Y)$ according
to their number of prime divisors as
\[
 B_{M,R}(Y)
 =
 \sum_{r=0}^R M^rA_r(Y)
 =
 1+\sum_{r=1}^R M^rA_r(Y).
\]
The initial term is $1$ because the only positive integer
with no prime divisors is $n=1$. Applying the bound \eqref{bound}
 gives
 \[
 B_{M,R}(Y)
 \ll_{M,R}
 1+\frac{Y}{\log Y}
 \sum_{r=1}^R M^r(1+\log\log Y)^{r-1}.
\]
For $Y\geq3$, we have $1+\log\log Y\geq1$. Hence
\[
 \sum_{r=1}^R M^r(1+\log\log Y)^{r-1}
 \leq
 \left(\sum_{r=1}^R M^r\right)
 (1+\log\log Y)^{R-1}.
\]
The finite sum depends only on $M$ and $R$. Moreover,
$Y/\log Y\geq1$, so the initial constant term satisfies
the same upper bound. This proves
\eqref{eq:colored-landau}. Finally, for every integer occurring in the defining
sum, $\omega(n)\leq R$ and therefore
$M^{\omega(n)}\leq M^R$. Since there are at most
$\lfloor Y\rfloor$ positive integers up to $Y$,
\[
 B_{M,R}(Y)
 \leq M^R\lfloor Y\rfloor
 \leq M^RY
 \ll_{M,R}Y,
\]
for $Y\geq 1$.
\end{proof}

In the applications in the next section we take $M=|I(G)|$ and
$R=R_G(t)$. 

\section{Counting nilpotent extensions of fixed genus}\label{sec:counting}

In this section, we prove Theorem~\ref{thm:main}. Recall that $G$ is a
fixed non-trivial finite nilpotent group, $\ell$ is the smallest prime
divisor of $|G|$,
\[
 I(G)=\{g\in G:\ord(g)=\ell\}\quad\text{and}\quad
 H(G)=I(G)\cup\{1\}.
\]
We also write
\[
 J(G):=G\setminus H(G)
      =\{g\in G:\ord(g)>\ell\}.
\]
For $h\in I(G)$ one has $e_h=e_0=a(G)^{-1}$, whereas
$a(G)e_g>1$ for every $g\in J(G)$.

For a real number $X\geq1$ and a positive integer $t$, set
\begin{equation}\label{eq:E-Gt-def}
 \cE_{G,t}(X)
 :=
 \left\{
 \psi\in\Epitop(\GQ,G):
 |\Disc(L_\psi)|\leq X,
 \ g_{L_\psi}=t
 \right\},
\end{equation}
where $L_\psi$ denotes the fixed field of $\ker(\psi)$.

\begin{lemma}\label{lem:epi-field-count}
For every $X\geq1$ and every positive integer $t$,
\begin{equation}\label{eq:epi-field-count}
 |\cE_{G,t}(X)|=|\Aut(G)|N_{G,t}(X).
\end{equation}
\end{lemma}

\begin{proof}
Apply Lemma~\ref{lem:aut-multiplicity} to the fields counted by $N_{G,t}(X)$.
\end{proof}

Fix $t$ and set
\begin{equation}\label{eq:R-fixed-section4}
 R:=R_G(t)=d_\ell(G)+\nu_\ell(t).
\end{equation}
By Corollary~\ref{cor:omega-fixed}, every $\psi\in\cE_{G,t}(X)$ has
\begin{equation}\label{eq:omega-min-section4}
 \omega_{\min,S_G}(L_\psi)\leq R.
\end{equation}
We translate this restriction into the Koymans--Pagano coordinates introduced in the previous section.

Let $v=(v_g)_{g\neq1}$ be the bundled primitive tuple corresponding to
$\psi$. Write $v_g=(v_g(1),v_g(2))$ as in
Section~\ref{sec:KP-local}. Since the first coordinates $v_g(1)$ take
values in fixed finite sets, the total number of possibilities for
$(v_g(1))_{g\neq1}$ is bounded in terms of $G$ only. They will
therefore be absorbed into the implied constants below.

For the squarefree coordinates define
\begin{equation}\label{eq:m-and-D-section4}
 m(v):=\prod_{h\in I(G)}v_h(2)\quad \text{and}\quad
 D(v):=\prod_{g\in J(G)}v_g(2)^{e_g}.
\end{equation}
The $v_g(2)$ are squarefree positive integers, and they are pairwise
coprime. Hence $m(v)$ is squarefree and
\begin{equation}\label{eq:omega-m-section4}
 \omega(m(v))
 =\sum_{h\in I(G)}\omega(v_h(2))
 =\omega_{\min,S_G}(L_\psi)
 \leq R.
\end{equation}
Moreover, Proposition~\ref{prop:KP-disc-general} gives
\begin{equation}\label{eq:disc-monomial-section4}
 m(v)^{e_0}D(v)
 =\operatorname{free}_{S_G}|\Disc(L_\psi)|
 \leq |\Disc(L_\psi)|
 \leq X.
\end{equation}
Since $e_0=a(G)^{-1}$, it follows that
\begin{equation}\label{eq:m-bound-section4}
 m(v)\leq
 \left(\frac{X}{D(v)}\right)^{a(G)}.
\end{equation}

For an upper bound, we enlarge the parameter set. For $g\in J(G)$,
replace the squarefree integer $v_g(2)$ by an arbitrary positive
integer $n_g$, and discard all coprimality, congruence and solvability
conditions involving these coordinates. Set
\begin{equation}\label{eq:Dn-section4}
 \mathbf n=(n_g)_{g\in J(G)},
\quad\text{and}\quad
 D(\mathbf n):=\prod_{g\in J(G)}n_g^{e_g}.
\end{equation}
In our arguments, for minimal coordinates we retain only the information that their
product $m$ is squarefree and that $\omega(m)\leq R$. Once such an
$m$ is fixed, each prime divisor of $m$ may be assigned to one of the
$|I(G)|$ coordinates indexed by $I(G)$. Thus the number of possible
ordered minimal-coordinate tuples with product $m$ is at most
$|I(G)|^{\omega(m)}$. Dropping the restrictions
$q\equiv1\pmod\ell$ and the Koymans--Pagano solvability conditions
can only increase the count.

Recall from Lemma~\ref{lem:colored-landau} that
\[
 B_{M,R}(Y)
 =
 \sum_{\substack{m\leq Y\\ \mu^2(m)=1\\ \omega(m)\leq R}}
 M^{\omega(m)}.
\]

\begin{lemma}\label{lem:master-section4}
One has
\begin{equation}\label{eq:master-section4}
 |\cE_{G,t}(X)|
 \ll_G
 \sum_{\substack{\mathbf n\\D(\mathbf n)\leq X}}
 B_{|I(G)|,R}
 \left(
 \left(\frac{X}{D(\mathbf n)}\right)^{a(G)}
 \right).
\end{equation}
\end{lemma}

\begin{proof}
By the Koymans--Pagano parametrization, it is enough to count the
corresponding solvable bundled primitive tuples. The first
coordinates have $O_G(1)$ possible values and may therefore be
absorbed into the implied constant.

Fix the nonminimal coordinates $\mathbf n$. The discriminant
inequality gives
\[
 m(v)^{e_0}D(\mathbf n)\leq X.
\]
Thus $D(\mathbf n)\leq X$ and, since $e_0=a(G)^{-1}$,
\[
 m(v)\leq
 \left(\frac{X}{D(\mathbf n)}\right)^{a(G)}.
\]
Furthermore, $m(v)$ is squarefree and $\omega(m(v))\leq R$.

For a fixed squarefree integer $m$, the ordered minimal-coordinate
tuples with product $m$ are obtained by assigning each prime divisor
of $m$ to one of the $|I(G)|$ minimal coordinates. Their number is
therefore at most $|I(G)|^{\omega(m)}$. Hence, for fixed
$\mathbf n$, the number of possible minimal-coordinate tuples is at
most
\[
 \sum_{\substack{
        m\leq (X/D(\mathbf n))^{a(G)}\\
        \mu^2(m)=1\\
        \omega(m)\leq R}}
 |I(G)|^{\omega(m)}
 =
 B_{|I(G)|,R}
 \left(
 \left(\frac{X}{D(\mathbf n)}\right)^{a(G)}
 \right).
\]Summing over $\mathbf n$
proves the result.
\end{proof}

Set $a=a(G)$. If $J(G)=\varnothing$, then
\eqref{eq:master-section4} reduces to
$B_{|I(G)|,R}(X^a)$, and Theorem~\ref{thm:main} follows immediately
from Lemma~\ref{lem:colored-landau}. Assume therefore that
$J(G)\neq\varnothing$, and set $r_G:=|J(G)|$. For $g\in J(G)$ define $\lambda_g:=ae_g$. By
Definition~\ref{def:eg}, one has $\lambda_g>1$. Set
\begin{equation}\label{eq:kappa-section4}
 \kappa_G
 :=
 \min_{g\in J(G)}
 \frac{ae_g-1}{r_Ge_g}>0.
\end{equation}
Choose an integer $B_G\geq1$ such that
\begin{equation}\label{eq:BG-section4}
 B_G\kappa_G\geq2,
\end{equation}
and set $Z_X:=(\log X)^{B_G}$. We take $X$ sufficiently large that
$Z_X<X$. For $D(\mathbf n)\leq X$, set
\[
 Y_{\mathbf n}:=
 \left(\frac{X}{D(\mathbf n)}\right)^a.
\]
Suppose first that $D(\mathbf n)\leq Z_X$. For $X$ sufficiently
large in terms of $G$,
\[
 \log Y_{\mathbf n}
 =a\log X-a\log D(\mathbf n)
 \geq \frac a2\log X.
\]
In particular, $Y_{\mathbf n}\geq3$ once $X$ is sufficiently large.
Also $1+\log\log Y_{\mathbf n}\ll_G1+\log\log X$. Hence
Lemma~\ref{lem:colored-landau} gives
\begin{equation}\label{eq:small-D-pointwise-section4}
 B_{|I(G)|,R}(Y_{\mathbf n})
 \ll_{G,t}
 \frac{X^a}{\log X}
 (1+\log\log X)^{R-1}
 D(\mathbf n)^{-a}.
\end{equation}
Summing over $D(\mathbf n)\leq Z_X$ and then dropping this restriction,
we obtain
\begin{align}
 &\sum_{D(\mathbf n)\leq Z_X}
 B_{|I(G)|,R}(Y_{\mathbf n}) \notag\\
 &\hspace{1cm}\ll_{G,t}
 \frac{X^a}{\log X}(1+\log\log X)^{R-1}
 \prod_{g\in J(G)}\zeta(ae_g).
 \label{eq:small-D-final-section4}
\end{align}
The product on the right is finite because $ae_g>1$ for every
$g\in J(G)$.

It remains to estimate the terms in \eqref{eq:master-section4} for
which $D(\mathbf n)>Z_X$. The crude estimate in
Lemma~\ref{lem:colored-landau} gives
\begin{equation}\label{eq:tail-first-section4}
 \sum_{Z_X<D(\mathbf n)\leq X}
 B_{|I(G)|,R}(Y_{\mathbf n})
 \ll_{G,t}
 X^a
 \sum_{D(\mathbf n)>Z_X}D(\mathbf n)^{-a}.
\end{equation}
If $D(\mathbf n)>Z_X$, then for some $g\in J(G)$,
$n_g^{e_g}>Z_X^{1/r_G}$. Consequently
\begin{align}
 \sum_{D(\mathbf n)>Z_X}D(\mathbf n)^{-a}
 &\leq
 \sum_{g\in J(G)}
 \left(\prod_{h\in J(G),\,h\neq g}\zeta(ae_h)\right)
 \sum_{n>Z_X^{1/(r_Ge_g)}}n^{-ae_g}.
 \label{eq:tail-union-section4}
\end{align}
For $N\geq1$, comparison with an integral gives
\[
 \sum_{n>N}n^{-ae_g}
 \leq N^{-ae_g}+\int_N^\infty u^{-ae_g}\,du
 \ll_G N^{1-ae_g},
\]
since $ae_g>1$. Taking $N=Z_X^{1/(r_Ge_g)}$, the exponent of $Z_X$
in this bound is
\[
 \frac{1-ae_g}{r_Ge_g}\leq-\kappa_G.
\]
It follows from \eqref{eq:kappa-section4} that
\begin{equation}\label{eq:tail-D-section4}
 \sum_{D(\mathbf n)>Z_X}D(\mathbf n)^{-a}
 \ll_G Z_X^{-\kappa_G}
 \ll_G(\log X)^{-2}.
\end{equation}
Combining \eqref{eq:tail-first-section4} and
\eqref{eq:tail-D-section4}, the tail is
$O_{G,t}(X^a(\log X)^{-2})$. Since $R\geq1$, this is smaller than the
right-hand side of Theorem~\ref{thm:main}.

\begin{proof}[Proof of Theorem~\ref{thm:main}]
If $J(G)=\varnothing$, apply Lemma~\ref{lem:colored-landau} directly
to \eqref{eq:master-section4}. If $J(G)\neq\varnothing$, combine
\eqref{eq:master-section4}, \eqref{eq:small-D-final-section4} and
\eqref{eq:tail-D-section4}. In both cases,
\[
 |\cE_{G,t}(X)|
 \ll_{G,t}
 \frac{X^{a(G)}}{\log X}
 (1+\log\log X)^{R_G(t)-1}.
\]
Theorem~\ref{thm:main} now follows from
Lemma~\ref{lem:epi-field-count}.
\end{proof}

\begin{corollary}\label{cor:density-proof}
For every fixed positive integer $t$,
\[
 \frac{N_{G,t}(X)}{N_G(X)}
 \ll_{G,t}
 \frac{(1+\log\log X)^{R_G(t)-1}}{\log X}.
\]
If, in addition, $I(G)\subseteq Z(G)$, then
\[
 \frac{N_{G,t}(X)}{N_G(X)}
 \ll_{G,t}
 \frac{(1+\log\log X)^{R_G(t)-1}}
 {(\log X)^{b(G,\Q)}}.
\]
\end{corollary}

\begin{proof}
The weak Malle lower bound of Kl\"uners--Malle gives
$N_G(X)\gg_GX^{a(G)}$; see \cite[Theorem~6.3]{KM}. Dividing the
estimate of Theorem~\ref{thm:main} by this lower bound proves the
first assertion. If $I(G)\subseteq Z(G)$, Koymans--Pagano prove
$N_G(X)\sim c_GX^{a(G)}(\log X)^{b(G,\Q)-1}$ with $c_G>0$; see
\cite[Theorem~1.1]{KP}. Division by this asymptotic gives the second
assertion.
\end{proof}

We now record the consequence of Corollary~\ref{cor:density-proof} for
Conjecture~1.4 of Dixit--Pasupulati \cite{DP}. Following their notation,
let
\[
 \cM(G,n,X)
 :=\{K/\Q:[K:\Q]=n,\ \Gal(\widetilde K/\Q)\simeq G,
                    \ |\Disc(K)|\leq X\},
\]
where $\widetilde K$ is the Galois closure of $K/\Q$ and fields are
counted up to $\Q$-isomorphism. Let $\cM'(G,n,X)$ be the subset
consisting of fields which admit a tower from $\Q$ in which every
successive extension is abelian.

\begin{corollary}\label{cor:DP}
Conjecture~1.4 of Dixit--Pasupulati holds for every nontrivial finite
nilpotent group in the regular representation. More precisely, for
every nontrivial finite nilpotent group $G$ and every fixed positive
integer $t$,
\[
 \frac{\#\{K\in\cM'(G,|G|,X):g_K=t\}}
 {\#\cM'(G,|G|,X)}
 \ll_{G,t}
 \frac{(1+\log\log X)^{R_G(t)-1}}{\log X}.
\]
In particular, this ratio tends to zero as $X\to\infty$.
\end{corollary}

\begin{proof}
If $K\in\cM(G,|G|,X)$ and $\widetilde K$ is its Galois closure, then
\[
 [K:\Q]=|G|=[\widetilde K:\Q],
\]
so $K=\widetilde K$. Thus the fields in $\cM(G,|G|,X)$ are precisely
the Galois $G$-extensions counted by $N_G(X)$. Since $G$ is
nilpotent, choose a subnormal series
\[
 1=G_m\triangleleft G_{m-1}\triangleleft\cdots\triangleleft G_0=G
\]
with cyclic prime-order quotients. The corresponding fixed fields
form a tower
\[
 \Q=K^{G_0}\subset K^{G_1}\subset\cdots\subset K^{G_m}=K,
\]
and $\Gal(K^{G_{i+1}}/K^{G_i})\simeq G_i/G_{i+1}$ is abelian
for every $i$. Thus every such extension belongs to
$\cM'(G,|G|,X)$. Hence
\[
 \cM'(G,|G|,X)=\cM(G,|G|,X),
\]
and the ratio in the statement is $N_{G,t}(X)/N_G(X)$. The result now
follows from Corollary~\ref{cor:density-proof}.
\end{proof}

\section{Bounds for genus moments}
\label{sec:weighted}
\label{sec:moments}

In this section, we prove Theorem~\ref{thm:moments-intro}. We first
obtain an upper bound for the moments for an arbitrary finite nilpotent
group. Under the centrality assumption, we construct a family of
extensions which gives a matching lower bound.

Throughout this section, $G$ is a nontrivial finite nilpotent group,
$\ell$ is the smallest prime divisor of $|G|$, and ${r}\geq0$ is fixed.
All estimates in $X$ are understood as $X\to\infty$. Set
\[
 a:=a(G),\qquad M:=|I(G)|,\qquad
 i:=i(G,\Q)=\frac{M}{\ell-1},\qquad u:=\ell^{r}.
\]
Recall that $i=b(G,\Q)$ when $I(G)\subseteq Z(G)$.

\begin{proposition}\label{prop:genus-upper-product}
Let $K/\Q$ be a finite Galois extension. Then
\begin{equation}\label{eq:genus-upper-product}
 g_K\leq\prod_{q\mid\Disc(K)}e_q(K/\Q).
\end{equation}
\end{proposition}

\begin{proof}
This follows from the narrow genus-number formula
\cite[Theorem~7]{Gurak}, applied over $\Q$.
\end{proof}

To estimate the moments, we compare the genus number with quantities
determined by ramification. The primes with inertia of order $\ell$
give a lower bound, while the product of the ramification indices gives
an upper bound.
For a Galois $G$-extension $L/\Q$, define
\begin{align}
 U_{r}(L)&:=\ell^{{r}\omega_{\min,S_G}(L)},
 \label{eq:Us-def}\\
 V_{r}(L)&:=\prod_{\substack{q\notin S_G\\q\text{ ramified in }L}}
                   e_q(L/\Q)^{r}.
 \label{eq:Vs-def}
\end{align}
If $v=(v_g)_{g\neq1}$ is its Koymans--Pagano tuple, then
Proposition~\ref{prop:bundled-inertia} gives
\begin{equation}\label{eq:Vs-coordinate-formula}
 V_{r}(L)=\prod_{g\neq1}\ord(g)^{{r}\omega(v_g(2))}.
\end{equation}
By Theorem~\ref{thm:genus-minimal} and
Proposition~\ref{prop:genus-upper-product},
\begin{equation}\label{eq:genus-weight-comparison}
 \ell^{-{r}d_\ell(G)}U_{r}(L)
 \leq g_L^{r}
 \leq |G|^{{r}|S_G|}V_{r}(L).
\end{equation}
For the upper bound, we used $e_p(L/\Q)\leq|G|$ at each
$p\in S_G$.
Thus, it suffices to obtain an upper bound for the sum of $V_{r}(L)$
and a lower bound for the sum of $U_{r}(L)$.

We shall use the following consequence of the Selberg--Delange method
to estimate the contribution from primes occurring in the minimal coordinates.

\begin{lemma}\label{lem:weighted-SD}
Let $T$ be a finite set of rational primes and let $z>0$. Put
\[
 B_{\ell,T,z}(Y):=
 \sum_{\substack{m\leq Y\\\mu^2(m)=1\\
 q\mid m\Rightarrow q\notin T,\ q\equiv1\; (\bmod\ell)}}
 z^{\omega(m)}.
\]
Then
\begin{equation}\label{eq:weighted-SD}
 B_{\ell,T,z}(Y)
 \asymp_{\ell,T,z}
 Y(\log Y)^{z/(\ell-1)-1}
 \qquad(Y\to\infty).
\end{equation}
In particular, for every $Y\geq1$,
\begin{equation}\label{eq:minimal-weighted-upper}
 B_{\ell,T,z}(Y)
 \ll_{\ell,T,z}
 Y(1+\log Y)^{z/(\ell-1)-1}.
\end{equation}
\end{lemma}

\begin{proof}
The associated Dirichlet series is
\[
 \prod_{\substack{q\notin T\\q\equiv1\; (\bmod\ell)}}
 (1+zq^{-w}),\qquad \Re(w)>1.
\]
Apply \cite[Theorem~4.1]{KP} to $\Q(\zeta_\ell)/\Q$, the
identity Frobenius class, and the excluded-prime set
$T\cup\{\ell\}$. The Selberg--Delange parameter is
$z/(\ell-1)$, and the leading constant is positive. This proves
\eqref{eq:weighted-SD}; enlarging the implied constant for bounded
$Y$ gives \eqref{eq:minimal-weighted-upper}.
\end{proof}

\begin{theorem}\label{thm:weighted-upper}
For every nontrivial finite nilpotent group $G$,
\begin{equation}\label{eq:Vs-upper-general}
 \sum_{\substack{L/\Q\;\mathrm{Galois}\\
 \Gal(L/\Q)\simeq G\\|\Disc(L)|\leq X}}V_{r}(L)
 \ll_{G,{r}}X^a(\log X)^{ui-1}.
\end{equation}
Consequently,
\begin{equation}\label{eq:genus-moment-upper-general}
 M_{G,{r}}(X)\ll_{G,{r}}X^a(\log X)^{ui-1}.
\end{equation}
\end{theorem}

\begin{proof}
We follow the counting argument of Section~\ref{sec:counting}, keeping
track of the weight $V_{r}(L)$.
Set $J:=J(G)=G\setminus(I(G)\cup\{1\})$. For the nonminimal
coordinates $\mathbf n=(n_g)_{g\in J}$, put
\[
 D(\mathbf n):=\prod_{g\in J}n_g^{e_g},\qquad
 c(\mathbf n):=\prod_{g\in J}\ord(g)^{{r}\omega(n_g)}.
\]
The first components $v_g(1)$ of the Koymans--Pagano tuple have $O_G(1)$ possible values. We retain
squarefreeness and pairwise coprimality among the $n_g$, and discard
the other restrictions on the nonminimal coordinates.

The product $m$ of the minimal coordinates is squarefree, supported
on primes $q\notin S_G$ with $q\equiv1\pmod\ell$, and satisfies
\[
 m^{1/a}D(\mathbf n)\leq X.
\]
Each prime dividing $m$ may be assigned to one of the $M$ coordinates
indexed by $I(G)$, and each assignment contributes the weight $u$. Discarding the solvability
conditions and the coprimality conditions between minimal and
nonminimal coordinates therefore gives
\begin{align*}
 &\sum_{\substack{L/\Q\;\mathrm{Galois}\\
 \Gal(L/\Q)\simeq G\\|\Disc(L)|\leq X}}V_{r}(L)\ll_G
 \sum_{\substack{\mathbf n\\D(\mathbf n)\leq X}}^{*}
 c(\mathbf n)
 B_{\ell,S_G,uM}\!\left((X/D(\mathbf n))^a\right),
\end{align*}
where the asterisk indicates squarefree, pairwise coprime
coordinates. The factor $|\Aut(G)|$ relating epimorphisms to fields
is absorbed into the implied constant.

Put $\kappa:=ui-1$. Since $M\geq\ell-1$ and $u\geq1$, we have
$\kappa\geq0$. If $D(\mathbf n)\leq X$ and
$Y=(X/D(\mathbf n))^a$, then
\[
 Y\geq1,\qquad
 (1+\log Y)^\kappa\ll_{G,{r}}(\log X)^\kappa
 \qquad(X\geq3).
\]
Applying Lemma~\ref{lem:weighted-SD} and extending the sum to all
squarefree, pairwise coprime tuples $\mathbf n$, with no restriction
on their prime divisors, we obtain
\begin{align*}
 \sum_{\substack{L/\Q\;\mathrm{Galois}\\
 \Gal(L/\Q)\simeq G\\|\Disc(L)|\leq X}}V_{r}(L)
 &\ll_{G,{r}}
 X^a(\log X)^\kappa
 \sum_{\mathbf n}^{*}\frac{c(\mathbf n)}{D(\mathbf n)^a}\\
 &=
 X^a(\log X)^\kappa
 \prod_p\left(1+\sum_{g\in J}
                   \frac{\ord(g)^{r}}{p^{ae_g}}\right).
\end{align*}
The Euler product converges because $ae_g>1$ for every $g\in J$.
If $J$ is empty, it is interpreted as $1$. This proves
\eqref{eq:Vs-upper-general}, and
\eqref{eq:genus-weight-comparison} gives
\eqref{eq:genus-moment-upper-general}.
\end{proof}

\begin{proposition}\label{prop:weighted-lower}
Suppose that $I(G)\subseteq Z(G)$. Then
\begin{equation}\label{eq:Us-lower}
 \sum_{\substack{L/\Q\;\mathrm{Galois}\\
 \Gal(L/\Q)\simeq G\\|\Disc(L)|\leq X}}U_{r}(L)
 \gg_{G,{r}}X^a(\log X)^{ub(G,\Q)-1}.
\end{equation}
\end{proposition}

\begin{proof}
For the lower bound, we fix one $G$-extension and twist its defining
epimorphism by characters taking values in the central subgroup
$H=I(G)\cup\{1\}$. These twists allow us to introduce ramification
of order $\ell$ at additional primes.

The set $H:=I(G)\cup\{1\}$ is a central elementary abelian
$\ell$-subgroup: for $h,h'\in H$, one has
$(hh')^\ell=h^\ell(h')^\ell=1$. Write $H$ additively, and let
$\iota:H\hookrightarrow G$ denote the inclusion.

Fix a $G$-extension $L_0/\Q$ once for each $G$, and let
$\psi_0:G_{\Q}\twoheadrightarrow G$ be an associated epimorphism.
Such an extension exists by Shafarevich's theorem; see
\cite[Theorem~6.1]{KM}. Let $T$ be the union of $S_G$ and the
finite primes ramified in $L_0$. Thus $T$ is fixed throughout the
construction. For each prime $q\notin T$ with
$q\equiv1\pmod\ell$, choose a nonzero character
\[
 \chi_q:G_{\Q}\longrightarrow\F_\ell
\]
corresponding to the degree-$\ell$ subfield of $\Q(\zeta_q)$.
It is ramified at $q$ and unramified at every other finite prime.

Let $m$ be squarefree and supported on these primes. For each
$q\mid m$, choose $h_q\in H\setminus\{0\}=I(G)$, and define
\[
 \xi:=\sum_{q\mid m}h_q\chi_q,
 \qquad
 \psi_\xi(\sigma):=\iota(\xi(\sigma))\psi_0(\sigma).
\]
Centrality of $H$ implies that $\psi_\xi$ is a homomorphism.
At every prime in $T$, its restriction to inertia agrees with that
of $\psi_0$. At a prime $q\mid m$, its inertia image is
$\langle h_q\rangle$, of order $\ell$.

We show that $\psi_\xi$ is surjective. Put
$J_\xi:=\operatorname{im}(\psi_\xi)$. Its image in $G/H$ is
surjective, so $J_\xi H=G$. Since $H$ is central, $J_\xi$ is
normal in $G$. The composite
\[
 G_{\Q}\xrightarrow{\psi_0}G\longrightarrow G/J_\xi
\]
is unramified at every finite prime. Indeed, at primes in $T$,
the inertia image of $\psi_0$ equals that of $\psi_\xi$ and is
contained in $J_\xi$; outside $T$, the map $\psi_0$ is unramified.
Minkowski's discriminant theorem implies that this quotient is
trivial. Thus $J_\xi=G$.

Let $L_\xi$ be the field defined by $\psi_\xi$. Outside $T$, it
is ramified precisely at the primes dividing $m$, each with tame
inertia of order $\ell$. The discriminant contribution of these
primes is therefore $m^{1/a}$. For $p\in T$, let
$E=(L_\xi)_w$ for a prime $w\mid p$. Since $L_\xi/\Q$ is Galois,
\[
 v_p\!\left(|\Disc(L_\xi)|\right)
 =
 \frac{|G|}{[E:\Q_p]}\,
 v_p\!\left(\Disc(E/\Q_p)\right).
\]
There are only finitely many extensions of $\Q_p$ of degree at
most $|G|$. Consequently, there is a constant $C_{G,T}\geq1$,
independent of $\xi$, such that
\[
 |\Disc(L_\xi)|\leq C_{G,T}m^{1/a}.
\]
Also, since $T\supseteq S_G$ and $u\geq1$,
\[
 U_{r}(L_\xi)\geq u^{\omega(m)}.
\]

We now count the epimorphisms obtained in this way. For $q\mid m$,
fix an inertia subgroup $I_q\subseteq G_{\Q}$. Since $\psi_0$ is
unramified at $q$, we have
\[
 \psi_\xi(\tau)=\iota\!\left(h_q\chi_q(\tau)\right)
 \qquad(\tau\in I_q).
\]
The restriction of $\chi_q$ to $I_q$ is surjective and has been fixed
in advance. Thus, different choices of $h_q$ give different restrictions
to inertia, even when they generate the same subgroup of $H$.
Different choices of $m$ give different ramified primes outside $T$.
Consequently, distinct choices of $m$ and $(h_q)_{q\mid m}$ give
distinct epimorphisms.
For each $m$ there are $M^{\omega(m)}$ such choices. Every field
arises from at most $|\Aut(G)|$ of these epimorphisms. Hence
\begin{align*}
 \sum_{\substack{L/\Q\;\mathrm{Galois}\\
 \Gal(L/\Q)\simeq G\\|\Disc(L)|\leq X}}U_{r}(L) &\geq\frac1{|\Aut(G)|}
B_{\ell,T,uM}\!\left((X/C_{G,T})^a\right)\\
 &\gg_{G,{r}}X^a(\log X)^{uM/(\ell-1)-1},
\end{align*}
by Lemma~\ref{lem:weighted-SD}. Since
$M/(\ell-1)=b(G,\Q)$, this proves \eqref{eq:Us-lower}.
\end{proof}

Combining the upper and lower bounds, we obtain the following result,
which proves Theorem~\ref{thm:moments-intro}.

\begin{theorem}\label{thm:genus-moments-central}
Suppose that $I(G)\subseteq Z(G)$. For every fixed real ${r}\geq0$,
\begin{equation}\label{eq:genus-moments-central}
 M_{G,{r}}(X)
 \asymp_{G,{r}}
 X^{a(G)}(\log X)^{\ell^{r}b(G,\Q)-1}.
\end{equation}
Moreover,
\begin{equation}\label{eq:average-moments-central}
 \frac{M_{G,{r}}(X)}{N_G(X)}
 \asymp_{G,{r}}
 (\log X)^{(\ell^{r}-1)b(G,\Q)}.
\end{equation}
\end{theorem}

\begin{proof}
By \eqref{eq:genus-weight-comparison}, the required bounds follow
from Theorem~\ref{thm:weighted-upper} and
Proposition~\ref{prop:weighted-lower}, since $i=b(G,\Q)$ under
the centrality assumption.
Taking ${r}=0$ gives
\[
 N_G(X)=M_{G,0}(X)
 \asymp_GX^{a(G)}(\log X)^{b(G,\Q)-1}.
\]
Dividing the two estimates proves
\eqref{eq:average-moments-central}.
\end{proof}

\newcommand{\etalchar}[1]{$^{#1}$}

\end{document}